\documentclass[11pt, twoside]{article}
\usepackage{amsmath}    
\usepackage{amssymb}   	
\usepackage{amsthm}     
\usepackage{graphicx}    

\usepackage{textcomp}    
\usepackage[T1]{fontenc} 
\usepackage{marvosym}  

\usepackage{authblk} 	
\usepackage[usenames]{xcolor} 
\usepackage{lastpage} 	

\usepackage[latin1]{inputenc} 
\usepackage{indentfirst}
\usepackage{amsfonts,amscd}	
\usepackage{empheq}
\usepackage{latexsym}		
\usepackage{graphicx,graphics,epsfig}								
\usepackage{graphpap}										
\usepackage{float}
\usepackage{xcolor}
\usepackage{makeidx}      
\usepackage{appendix}	  
\usepackage{multicol}										
\usepackage{color}

\usepackage[normalem]{ulem}
\usepackage{changepage}
\usepackage{hyperref} 
\hypersetup{
	colorlinks=true,
	linkcolor=blue,
   citecolor=blue,
   urlcolor=blue
}

\usepackage[dvips=false,pdftex=false,vtex=false]{geometry}
\usepackage{dsfont} 

\numberwithin{equation}{section}
\numberwithin{table}{section}
\numberwithin{figure}{section}

\newcommand{\orcid}[1]{\href{https://orcid.org/#1}{\includegraphics{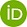}}} 

\newtheorem{theorem}{Theorem}[section]
\newtheorem{corollary}[theorem]{Corollary}
\newtheorem{lemma}[theorem]{Lemma}

\theoremstyle{definition}
\newtheorem{definition}[theorem]{Definition}
\theoremstyle{definition}

\theoremstyle{definition}

\title{A survey on obstacle-type problems for fourth order elliptic operators}

\author[1]{\textbf{Donatella Danielli}\footnote{e-mail: ddanielli@asu.edu}}
\author[2]{\textbf{Alaa Haj Ali}\footnote{e-mail: ahajali1@asu.edu}}

\affil[1]{Arizona State University, School of Mathematical and Statistical Sciences, Tempe, AZ,  85287-1804, USA. }

\affil[2]{Arizona State University, School of Mathematical and Statistical Sciences, Tempe, AZ,  85287-1804, USA. }

\newcommand\shorttitle{}

\makeatletter
\def\volume#1{\def\@volume{#1}}

\def\shortauthor#1{\def\@shortauthor{#1}}

\def\shorttitle#1{\def\@shorttitle{#1}}

\makeatother

\date{2021} 
\volume{\#} 

\makeatletter
\let\thetitle\@title
\let\thedate\@date
\let\thevol\@volume
\let\theauthor\@author
\let\theshortauthor\@shortauthor
\let\theshorttitle\@shorttitle
\makeatother

\makeatletter
\renewcommand{\maketitle}{\bgroup\setlength{\parindent}{0pt}

\parindent=1em
\renewcommand{\thefootnote}

\phantom{ }

\vspace{-1.5cm} \noindent\includegraphics[scale=.35]{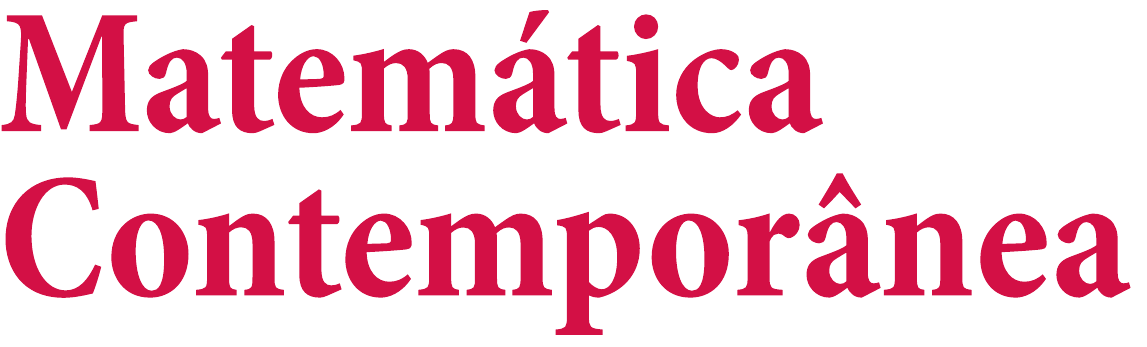}
 \hfill
 \includegraphics[scale=0.07]{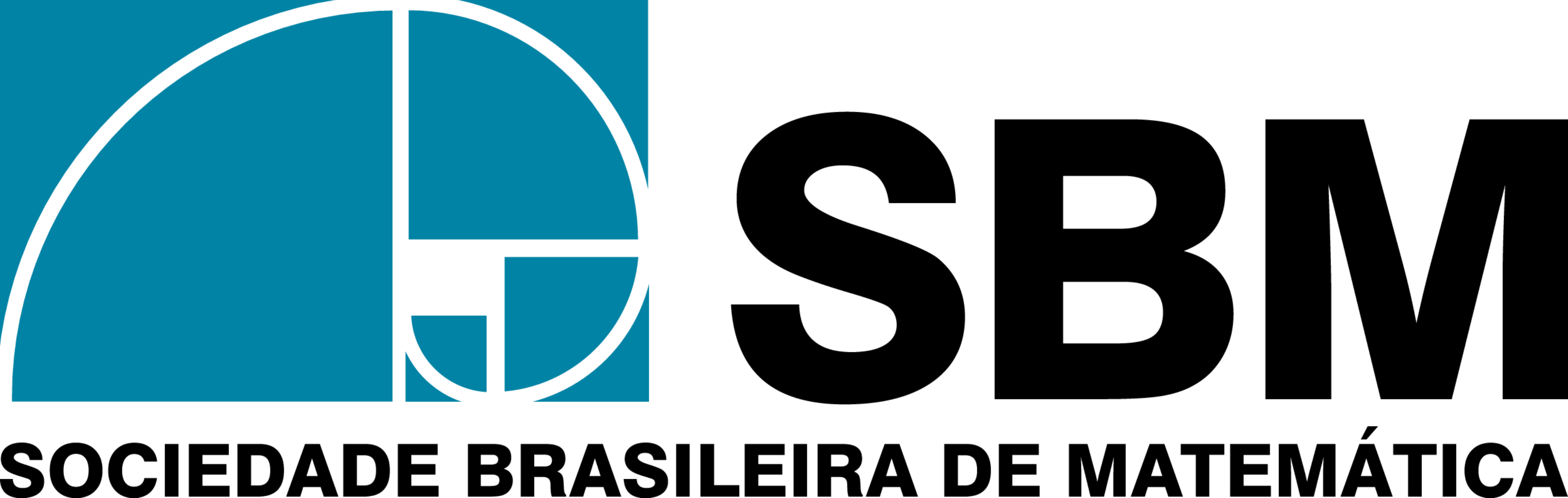}

\vspace*{-.2cm}

\noindent{\scriptsize Vol. \thevol,\  \pageref{FirstPage}--\pageref{LastPage}} \hfill \copyright \thedate \qquad \qquad

\vspace*{-.2cm}

\noindent{\scriptsize\url{http://doi.org/will\_be\_added\_later}}
\vspace{1truecm}
\begin{center}{\vbox{\titlefont\@title}}\end{center}
\vspace{0.5truecm}
\begin{center}{\@author} \end{center}

\egroup
}

\renewcommand{\thefootnote}{\fnsymbol{footnote}}
\renewcommand{\@fnsymbol}[1]{%
    \ifcase#1 \or {\,\Letter\!} \or\textasteriskcentered\or \textasteriskcentered\textasteriskcentered
    \else\@ctrerr\fi}
\makeatother

\newcommand*{\titlefont}{\fontsize{18}{21.6}\selectfont\textbf}

\makeatletter
\renewcommand\@author{\ifx\AB@affillist\AB@empty\AB@author\else
      \ifnum\value{affil}>\value{Maxaffil}\def\rlap##1{##1}%
    \AB@authlist\\[\affilsep]\vbox{\AB@affillist}
    \else  \AB@authors\fi\fi}
\makeatother

\newcommand\R{\mathds{R}}

\begin{document}

\maketitle
\thispagestyle{plain}
\renewcommand{\thefootnote}{\arabic{footnote}}   
\setcounter{footnote}{0}     
\setcounter{page}{1} 
\label{FirstPage}	

\pagestyle{myheadings} \markboth{\hfil  $\hspace{1.5cm}$  D. Danielli and A. Haj Ali
\hfil $\hspace{3cm}$ } {\hfil$\hspace{1.5cm}$
{ Obstacle problems for fourth order operators}
\hfil}

\begin{center}
\noindent
\begin{minipage}{0.85\textwidth}\parindent=15.5pt

\smallskip

{\small{
\noindent {\bf Abstract.} In this article we give a brief overview of some known results in the theory of obstacle-type problems associated with a class of fourth-order elliptic operators, and we highlight our recent work with collaborators in this direction. Obstacle-type problems governed by operators of fourth order naturally arise in the linearized Kirchhoff-Love theory for
plate bending phenomena. Moreover, as first observed by Yang in \cite{Y13},   boundary obstacle-type problems associated with the weighted bi-Laplace operator can be seen as extension problems, in the spirit of the one introduced by Caffarelli-Silvestre, for the fractional Laplacian $(-\Delta)^s$ in the case  $1<s<2$.
In our recent work, we investigate some problems of this type, where we are concerned with the well-posedness of the problem, the regularity of solutions, and the structure of the free boundary. In our approach, we combine classical techniques from potential theory and the calculus of variations
with more modern methods, such as the localization of the operator and monotonicity formulas.}
\smallskip

\smallskip

}

\end{minipage}
\end{center}



\section{Introduction} 

The aim of this survey is to introduce the interested reader to the theory of obstacle-type problems involving a fourth-order operator.

Obstacle-type problems for second order uniformly elliptic operators have attracted  great interest over the years, and they continue to be explored. Such problems have a wide range of applications in elasticity, biology, engineering, and mathematical finance, just to name a few.  The quintessential prototype is the classical obstacle problem, which consists in finding the equilibrium position of an elastic membrane whose boundary is held fixed, and which is constrained to lie above a given obstacle. Mathematically, given a domain $\Omega$, a function $\psi$ (which represents the obstacle), and a function $g$ (given boundary data) satisfying the compatibility condition $g \geq \psi$ on $\partial \Omega$, one seek to find the membrane's vertical displacement function which minimizes the Dirichlet energy
\begin{equation}\label{fnct_lap}
J[w]=\int_{\Omega} |\nabla w|^2,
\end{equation}
over all functions $w$ satisfying $w=g$ on $\partial \Omega$, and constrained to remain above the obstacle $\psi$ in $\Omega$.

Another example is the  thin obstacle problem, also known as the \textit{Signorini problem}, first proposed and studied by Signorini \cite{Sig59} in connection with linear elasticity. In mathematical terms, given a domain $\Omega$, and a lower dimensional subset $\Gamma$ of $\Omega$, the functional (\ref{fnct_lap}) is minimized over functions $w$ satisfying $w=g$ on $\partial \Omega$, and constrained in this case to remain above a given obstacle only on the thin set $\Gamma$. The same mathematical model also appears in the study of semipermeable media and boundary heat control. Moreover, it can be seen a local formulation of an  obstacle problem for the fractional Laplacian which has important implications in financial mathematics.

Obstacle-type problems feature an a priori unknown interface which is defined as the topological boundary (in the domain in consideration for global obstacle problems, and in a lower dimensional subdomain for thin obstacle problems) of the coincidence set between the solution and the given obstacle. The main goal in studying such problems is to understand regularity properties of the solution, as well as the structure of the free interface.  Obstacle-type problems for second order uniformly elliptic operators have been extensively studied over the years. These efforts have been spearheaded by the innovative  blow-up techniques introduced by Caffarelli in his seminal paper \cite{Caf77}, and by the discovery of several powerful monotonicity formulas which provided the necessary tools to establish the optimal regularity of solutions, as well as the regularity and the structure of the free boundary. Regarding the classical obstacle problem, we refer the reader to \cite{Fre72} \cite{Caf77}, \cite{Caf98}, \cite{KN77}, \cite{Mon03}, \cite{Wei99}, \cite{FS19}, whereas for thin obstacle problems, we refer to \cite{ACS08}, \cite{AC04}, \cite{CSS08}, \cite{CS07}, \cite{GP09}, the book \cite{PSU12} and the survey \cite{DS17}, as well as the references therein. These lists are, of course, far from being exhaustive.\\
\indent
On the other hand, much less work has been done on free boundary problems associated with the bi-Laplace operator. Regarding obstacle-type free boundary problems, the regularity of the solution has been studied in \cite{Fre71}, \cite{Fre73}, \cite{CF79} for the global obstacle case; and in \cite{Sc84}, and \cite{Sc86} for the thin obstacle case. Moreover, in \cite{CF79} the authors considered the structure of the free boundary of a solution, and  studied its regularity  in the two dimensional case. Also, in \cite{Ale19}, Aleksanyan investigated the regularity of the free boundary for the global obstacle problem under the assumptions that the solution is almost one-dimensional, and that the non-coincidence set is a non-tangentially accessible domain (NTA).  A related model is  the singularly perturbed bi-Laplace equation $ \Delta^2 u^{\epsilon}=-\beta_{\epsilon}(u^{\epsilon})$, which can be thought of  as the biharmonic counterpart of the classical combustion problem. This problem was investigated by Dipierro, Karakhanyan, and Valdinoci in \cite{DKV19}, where they proved the convergence as $\epsilon \to 0$ to a free boundary problem driven by the bi-Laplacian. Moreover, they derived a monotonicity formula in the plane, and used it to establish the quadratic behavior of solutions near the zero level set. The limiting problem (as $\epsilon \to 0$) was  investigated by the same authors in \cite{DKV20}, where
the authors prove that  at   free boundary points where $\nabla u(x)=0$, the free boundary   is either well approximated by zero-level sets of quadratic homogeneous polynomials, or $u$ has quadratic growth. \\
\indent
Many important questions related to the optimal regularity of a solution of an obstacle-type problem, and to the structure of its free boundary and the characterization of  free boundary points remain open. The main objective of this article is to give a brief survey of the related results available in the literature, with some focus on our recent contribution to this area of research.

\section{Obstacle problems for plates}\label{applications}
In the two dimensional case, problems associated with the bi-Laplace operator naturally arise  in elastic, homogenous and isotropic plate displacement phenomena in the event of relatively small displacements. In such cases, a three dimensional plate with given volume $\mathcal{V}$ and with small thickness $h$ is under a normal load force $f$, which is sufficiently weak so that induced displacements $(u_1,u_2,u_3)$ satisfy $u_1(x_1,x_2,0)=u_2(x_1,x_2,0)=0$. In this approximate theory, we assume that our plate is two dimensional occupying a region $\Omega$ in the plane $\{x_3=0\}$. Hence, we can now simplify the notation and write $u(x_1,x_2)$ for  $u_3(x_1,x_2,0)$. In the case of small deformations, the displacement $u$ is governed by the equation
\begin{equation*}
D\Delta^2 u =f.
\end{equation*}
Here,  $D=\frac{E h^3}{12(1-\delta^2)}$ is known as the \textit{stiffness coefficient of the plate} or \textit{modulus of flexural rigidity}, where $E$ is the Young's modulus, and $\delta$ is the Poisson's ratio.

Next, we describe some  possible boundary conditions, along with their physical interpretations. Our writing is based on Chapter 4 of the book \cite{DL76} by Duvaut and Lions (see also Chapter 2 of the book \cite{LL59} by Landau and Lifshitz). The boundary conditions  will involve either the displacement $u$ itself, or the given forces. For a point $(x_1,x_2)$ on the boundary  $\partial \Omega$ of our flat plate, we denote by $\nu=(\nu_1, \nu_2)$ the exterior unit normal to $\partial \Omega$, and by $\tau$ the unit tangent obtained from $\nu$ through a rotation by $\pi/2$. There are two components of the line density of moments acting on $\partial \Omega$:

\begin{equation}\label{M1}
M_1=  M_{x_1 x_2} \ \nu_1 + M_{x_2} \  \nu_2
\end{equation}
and
\begin{equation}\label{M2}
M_2=-M_{x_1}  \ \nu_1-  M_{x_1 x_2} \  \nu_2 .
\end{equation}
Here
\begin{equation*}\label{M_x}
\begin{split}
&M_{x_1}:=D \left(\frac{\partial^2 u}{\partial x_1^2} + \delta \frac{\partial^2 u}{\partial x_2^2}\right) \text { and } M_{x_2}:=D\left(\frac{\partial^2 u}{\partial x_2^2} + \delta \frac{\partial^2 u}{\partial x_1^2}\right) \\
&\text { are the bending moments, and }\\
&M_{x_1 x_2}:=D(1-\delta)\frac{\partial^2 u}{\partial x_1 \partial x_2} \text { is the twisting moment. }\\
\end{split}
\end{equation*}
One can naturally form the following normal and tangential components of this density of moments:
\begin{equation}\label{M_nu}
M_{\nu}:=M_1 \nu_1 + M_2 \nu_2,
\end{equation}

\begin{equation}\label{M_tau}
M_{\tau}:=M_2 \nu_1 -M
_2 \nu2,
\end{equation}
Moreover, suppose $\Omega_1$ is a portion of $\Omega$ with boundary $\Gamma_1$, and let $\mathcal{V}_1:=\Omega_1 \times (-h/2, h/2)$ be the volume  corresponding to $\Omega_1$. The portion $\mathcal{V} \setminus \mathcal{V}_1$ exerts a system of forces on $\mathcal{V}$ which, in equilibrium, together with the density $f$ of forces on $\Gamma_1$, form a system statically equivalent to zero. On the other hand, the forces exerted by $\mathcal{V} \setminus \mathcal{V}_1$ on $\mathcal{V}_1$ act across the interface $S_1 =\Gamma_1 \times (-h/2, h/2)$.

The formulation of equations in \cite{DL76} suggests that certain well posed problems would have  boundary conditions on $\partial \Omega$ involving
 either $F$ or $u$, and either $M_{\tau}$ or $\frac{\partial u}{\partial \nu}$.
Here,
\begin{equation}\label{F}
F=R-\frac{\partial M_{\nu}}{\partial \tau},
\end{equation}
where $R$ depends on the restriction to $\Gamma_1$ of the forces acting across the interface $S_1$ as described above.
For example, in the case of a clamped plate, we prescribe
\begin{equation}\label{ex1}
u = \frac{\partial u}{\partial \nu}=0 \text { on }\partial \Omega.
\end{equation}
This says that  the plate is clamped on its boundary, and it remains horizontal there.
In another setting, we would have a plate with unrestricted boundary on a portion $\Gamma_1 \subset \partial \Omega$ of its edge, and with fixed boundary on the remaining  of $\partial \Omega$. One then has
\begin{equation}\label{ex2}
F=M_{\tau}=0 \text { on } \Gamma_1, \text  { and } u=M_{\tau}=0 \text { on } \partial \Omega \setminus \Gamma_1 .
\end{equation}
\subsection{Unilateral displacements}

\subsubsection{Unilateral displacements of the points of $\Omega$}
This is a global obstacle problem, which can be stated as
\begin{equation}\label{thick_application}
u\geq 0, \ \Delta^2 u \geq f \text { in }\Omega,  \text { and } \Delta^2 u =f \text{ in } \{u>0\},
\end{equation}
with classical  boundary conditions on $\partial \Omega$ of the form (\ref{ex1}) or (\ref{ex2}). Here $f$ is the load force exerted on the plate.
\subsubsection{Unilateral displacements of the points of $\partial \Omega$}
This is the thin obstacle-type problem
\begin{equation}\label{thin_application}
\Delta^2 u =f \text { in } \Omega, \ \  F \geq 0, u \geq 0, \text { and }  u \  F=0 \text { on } \Gamma_1,
\end{equation}
where $\Gamma_1$ is a part of the boundary $\partial \Omega$.
The conditions on $M_{\tau}$  or $\frac{\partial u}{\partial \nu}$ on both of $\Gamma_1$ and $\partial \Omega \setminus \Gamma_1$ are of classical type as in (\ref{ex1}) or (\ref{ex2}). Here $f$ is the load force, and $M_{\tau}$  and $F$ are as defined in (\ref{M_tau}) and (\ref{F}) respectively.

\section{The classical obstacle problem for the bi-Laplacian}\label{thick_section}
For a smooth bounded domain $\Omega\subset\mathbb{R}^n$, and given functions $g,\psi\in C^2(\overline{\Omega})$ satisfying the compatibility condition $\psi \leq g$ on $\partial \Omega$ in some trace sense, consider the problem of minimizing the energy
\begin{equation}\label{fnct_classical}
\int_{\Omega} (\Delta w)^2 dx
\end{equation}
over the set
\begin{equation}\label{admissible_classical}
\mathcal{K}:=\{w \in H^2(\Omega) :  w-g \in H_0^2(\Omega) \text {, and } w \geq \psi \text{ a.e.  in } \Omega\}.
\end{equation}
Following a variational approach, one can show that the minimizer of (\ref{fnct_classical}) over the set (\ref{admissible_classical}) satisfies
\begin{equation}\label{prb_classical}
\left\{
\begin{array}{ll}
u \geq \psi & \text { in } \Omega\\
\Delta^2 u=0 & \text { in } \{u>\psi\} \\
\Delta^2 u \geq 0 & \text { in } \Omega.
\end{array}\right.
\end{equation}
In other words, the solution $u$ stays above the obstacle $\psi$, it is biharmonic whenever it does not touch the obstacle, and it is super-biharmonic in the whole domain. This splits the domain $\Omega$ into two regions: the \textit{non-contact set} $\{u>\psi\}$ in which the solution $u$ is biharmonic, and the \textit{contact set} $\{u=\psi\}$. The interface in between the two regions is called the \textit{free boundary} which depends on the solution, and hence it is a priori unknown.

\subsection{Regularity of solutions}
In \cite{Fre71}, Frehse established the $H^3_{loc}$-regularity of the solution $u$ to problem \eqref{prb_classical}. Then, in \cite{Fre73}, he proved that the solution is in  $C^{1,1}_{loc}(\Omega)$.

In \cite{CF79}, Caffarelli and Friedman took a different approach, which hinges on techniques from potential theory, to prove the $C^{1,1}_{loc}$-regularity of a solution. Using a variational argument, the authors showed that $\mu:=\Delta^2 u$ defines a measure which is finite on compact subsets of the domain $\Omega$. Also, assuming that $\phi<0$ on $\partial \Omega$, the authors proved that the total mass $\mu(\Omega)$ is finite. To reach this result, the authors constructed some approximating problems with solutions $u^{\epsilon}$ converging, up to a subsequence, to the solution $u$, and then they proved that the finite measures induced by $u^{\epsilon}$ converge to $\mu$ in some weak sense. Moreover, in the case where $n \leq 4$, and still under the hypothesis $\phi<0$ on $\partial \Omega$, the authors proved that a solution is smooth up to the boundary. In addition, they showed that, in the two-dimensional case,  $u\in C^2(\Omega)$. Finally,  a counter-example is produced to show that an a-priori estimate on the modulus of continuity of $D^2 u$ in a compact subset of $\Omega$ cannot hold.

We now present a few salient details on the proof in \cite{CF79} of the $C^{1,1}_{loc}$-regularity of a solution $u$, since we will return to this circle of ideas later on. By considering the mollifiers of $u$ as smooth approximation, and  applying  Green's formula to these mollifiers, the authors show that $\Delta u$ is equal a.e. to an upper semi-continuous function $w$, and they derive an estimate for $w$ on the contact set of the solution $u$. More precisely, they reach the following results.
\begin{lemma}\textup{(\cite[Lemma 2.1]{CF79})}\label{L1}
There exists an upper semicontinuous function $w$ such that $w=\Delta u$ a.e. in $\Omega$. Moreover, for every $x^0 \in \Omega$, and for any sequence of balls $B_{\rho}(x^0)$
\begin{equation*}
\frac{1}{|B_{\rho}(x^0)|}\int_{B_{\rho}(x^0)} w dx \searrow w(x^0) \text { as $\rho \searrow 0$}.
\end{equation*}
\end{lemma}

\begin{theorem}\label{L2*}\textup{(\cite[Theorem 2.2]{CF79})}
For any point $x^0$ on the support of the measure $\mu$, we have
\begin{equation}\label{L2}
 w(x^0) \geq \Delta \phi(x^0).
\end{equation}
\end{theorem}
In order to apply techniques from potential theory, the authors consider the fundamental solution of the Laplace equation for $n \geq 3$
\begin{equation}\label{fundamental_classical}
V(x):=\gamma |x|^{2-n}, \ \   (\gamma>0),
\end{equation}
and they derive the following relation (see the proof of Theorem 3.1 in \cite{CF79}), which expresses the function $w$ in terms of the convolution $V$ with $\mu_{\rho}$, the restriction of the measure $\mu$ to $B_{\rho}(x^0)$,
\begin{equation}\label{relation_classical}
w = -V* \mu_{\rho} + \delta_{\rho}(x),\text { for all $x \in B_{\frac{\rho}{2}}(x^0)$}.
\end{equation}
Here, $\delta_{\rho}(x)$ is a bounded function on $B_{\frac{\rho}{2}}(x^0)$.
With the equation (\ref{relation_classical}) in hand, to obtain a bound for $w$ (and hence for $\Delta u$), it suffices to obtain a bound for the convolution $\hat{V}:=V* \mu_{\rho}$. The authors observed that $\hat{V} \geq 0$, and combined the results of Lemma \ref{L2*} and  (\ref{relation_classical}) to reach an upper bound for $\hat{V}$ on the support of the measure $\mu_{\frac{\rho}{2}}$.
At this point, the authors point out that $\hat{V}$ superharmonic, so that they apply the following maximum principle to reach an upper bound for $\hat{V}$ in the whole domain $\Omega$.

\begin{theorem}\label{MP_classical} \textup{(\cite{Cal67} and \cite{Lan72})}
Let $Z$ be a superharmonic function in $\mathbb{R}^n$ for which the measure $\nu=-\Delta Z$ is supported on a bounded set $S$. If $Z \leq M$ on $S$, then $Z \leq M$ in all of $\mathbb{R}^n$.
\end{theorem}

\subsection{Regularity of the free boundary}
In \cite{CF79}, the authors carry some analysis of the regularity of the free boundary. Under the assumption $\Delta^2 \phi >0$, the authors prove that the non-contact set is connected. Moreover, in the two-dimensional case and under the assumption $\Delta u(x^0)>\Delta \phi(x^0)$, the authors show that the free boundary in a neighborhood of $x^0$ is contained in a continuously differentiable curve.

In \cite{Ale19}, Aleksanyan studied the regularity of the free boundary for the  biharmonic obstacle problem under the assumptions that the obstacle is zero,   the solution is almost one-dimensional, and  the non-coincidence set is a non-tangentially accessible domain (NTA). The author derives the $C^{1,\alpha}$-regularity of the free boundary following a linearization argument introduced by Andersson in \cite{And20}.

Given a free boundary point, which without loss of generality we suppose is the origin, and for a ball $B$ around the origin, the author defines the class $\mathcal{B}(\epsilon)$ as, among other requirements,  the set of almost one-dimensional functions in the sense $\|\nabla^{'} u\|_{H^2(B)}\leq \epsilon$ where $x^{'}=x-x_ne_n$ and  $\nabla^{'}=\nabla_x-u_{x_n} e_n$. The author shows that, for $\epsilon$ small enough, there exists a rescaling function
 \begin{equation*}
u_{r}(x)=\frac{u(rx)}{r^{3}} \  \ (x \in B),
\end{equation*}
and  a constant $0< \gamma<1$ such that
\begin{equation*}
\|\nabla^{'} u_{r}\|_{H^2(B)}\leq \gamma \|\nabla^{'} u\|_{H^2(B)}\leq \gamma \  \epsilon.
\end{equation*}
in a special coordinate system depending on $u$. Iterating the above inequality and changing the  coordinate system accordingly, the author establishes the existence of a measure theoretic normal $\eta_0$ to the free boundary $\Gamma_u$ of the solution at the origin.
For a general free boundary point $x$, letting $\eta_x$ be the measure theoretic normal to $\Gamma_u$ at $x$, the author shows that $\eta_x$ is a H\"older-continuous function on $\Gamma_u$.

Moreover, the author uses the above regularity of the free boundary result to improve the regularity of the solution to
 $C^{3,\alpha}_{loc}$, and therefore $C^{2,1}$. The proof is based on applications of the boundary regularity theory for harmonic functions, Calderon-Zygmund estimates, and a Sobolev embedding theorem. Finally, the author gives an example to show that, in general - that is, without the assumptions that the solution is almost one-dimensional and the non-coincidence set is an NTA domain - $C^{2,1/2}$ is the best regularity a solution may achieve in dimension $n \geq 2$.

\section{The thin obstacle problem for the poly-Laplacian}\label{thin_section}
Let $m \geq 2$ be an integer, and suppose $\Omega$ is a smooth domain in  $\mathbb{R}^{n+1}$ ($n \geq 1$, $2m \leq n+3$). Also, let  $f \in L^p(\Omega)$, with $p>\frac{n+1}{2m-1}$, and $g \in C^2(\overline{\Omega})$. Given an obstacle function $\Psi \in C^2(\overline{\Omega})$, the  \emph{classical obstacle problem} for the Laplacian of order $m$ consists in finding a function $u$ solution of the variational inequality
\begin{equation}\label{m_varineq_thick}
<(-\Delta)^m u -f, w-u> \geq 0
\end{equation}
for all $w$ satisfying $w-g \in H_0^m(\Omega)$ and $w \geq \Psi\text { on } \Omega$.

If $\Omega^{'}\subset \Omega$ is a $k$-dimensional $C^2$-surface ($k \geq n-2m+4$) such that $\partial \Omega^{'}\subset \partial \Omega$, given a \emph{thin} obstacle function $\psi \in C^2(\overline{\Omega'})$, the \emph{thin obstacle problem} for the Laplacian of order $m$ is governed by the variational inequality
\begin{equation}\label{m_varineq_thin}
<(-\Delta)^m u -f, w-u> \geq 0
\end{equation}
for all $w$ satisfying $w-g \in H_0^m(\Omega)$ and $w \geq \psi\text { on } \Omega'$. For $m=2$ and $f=0$, problem \ref{m_varineq_thick} reduces to the classical problem for the bi-Laplacian presented in Section \ref{thick_section}.

 In \cite{Sc84},  Schild considered both  problems (\ref{m_varineq_thick}) and (\ref{m_varineq_thin}) for a general dimension $n \geq 1$, and a general order $m\geq 2$ satisfying $2m \leq n+3$. His  approach  hinges on techniques from potential theory. First, he established the $C^{1,1}_{loc}$-regularity of a solution $u$ of the problem under consideration. Observing that $\mu:=(-\Delta)^m u -f$ defines a non-negative measure on the domain, the author considered the restrictions $\mu_{\rho}$ of the measure $\mu$ on  balls $B_{\rho}$, and wrote both the  solution and its Laplacian as a sum of a potential and a regular function. Then, to obtain estimates for the potential (and therefore for the solution and its Laplacian), Schild used tools from potential theory such as a continuity principle and a maximum principle, in addition to some arguments involving smooth approximations of functions and Green's formula. Moreover, a crucial step was observing that the mixed second order partial derivatives of the fundamental solution  of the bi-Laplace equation are controlled by the Laplacian of such solution. This opens the way for deriving estimates on the mixed second order partial derivatives of the solution to the problem in consideration. Moreover, adapting a classical result from potential theory which provides a relation between the $L^2$-norm of a Riesz potential of a measure and a Riesz energy of the same measure, the author was able to use the  $C^{1,1}_{loc}$-regularity of a solution to establish its $W_{loc}^{3,2}$-regularity. Finally, the author gave a counter-example for higher regularity for the two dimensional case to show that,  in general, one cannot expect the solution to be in $W_{loc}^{3,\infty}(\Omega)$. More details on this potential theory based approach will appear in Section \ref{paper1}, where we present our work in \cite{DHP}.

\section{The fractional Laplacian of order $1<s<2$}\label{prelim}

In this section, we aim to present the connection between some fourth order linear operators, and the fractional Laplacian $(-\Delta)^s$ for $1<s<2$. We start with recalling definitions of the related function spaces and the associated norms.

On the Schwartz space of smooth functions which rapidly decay at $\infty$,
\begin{equation}\label{S}
\mathcal{S}:=\left\{v \in C^{\infty}(\mathbb{R}^n) \  \text{ s.t. } \sup_{ \mathbb{R}^n} |x|^{m} v^{(k)}(x)<+\infty  \text { for all $m\geq 0$, $k\geq 0$}\right\},
\end{equation}
one can define the fractional Laplacian $(-\Delta)^s $ as a singular integral
\begin{equation}\label{singular_frac}
\begin{split}
(-\Delta)^s v(x) &=C(n,s) \textit { P.V.} \int_{\mathbb{R}^n} \frac{\displaystyle \nabla v(x)-\nabla v(z)}{\displaystyle |x-z|^{n-2+2s}} dz\\
&:=C(n,s) \lim_{r \to 0} \int_{\mathbb{R}^n \setminus B_r(x)} \frac{\nabla v(x)-\nabla v(z)}{\displaystyle |x-z|^{n-2+2s}} dz\\
\end{split}
\end{equation}
for all $x \in \mathbb{R}^n$.
Here,\textit{ P.V.} is an abbreviation of \textit{``principal value''} as defined in the latter equation. The definition (\ref{singular_frac}) is equivalent to the following one in terms of a Fourier multiplier
\begin{equation}\label{Fourier_frac}
\widehat{(-\Delta)^s v}= |\xi|^{2 s} \widehat{v}(\xi)   \text {\ \  \ \ for all $v \in \mathcal{S}$,}
\end{equation}
where $\hat{f}$ denotes the Fourier transform of $f$.   We introduce the  Gagliardo-Slobodeckii semi-norm as
\begin{equation}\label{semi-norm}
\begin{split}
[v]_{\dot H^s}&:=\left(\int_{\R^n} \int_{\R^n} \frac{|\nabla v(x)-\nabla v(z)|^2}{|x-z|^{n-2+2s}} dx dz\right)^{1/2}\\
&=\left(k(n,s)\int_{\R^n}|\xi|^{2s}|\hat v(\xi)|^2d\xi\right)^{1/2},
\end{split}
\end{equation}
where  $k(n,s)$ is a positive constant depending on the dimension $n$ and the exponent $s$ only.

We define the space $\dot{H}^s(\mathbb{R}^n)$ as the completion of the $C_0^{\infty}(\mathbb{R}^n)$ with respect to the semi-norm (\ref{semi-norm}).
We observe that
\begin{equation}\label{inner_equivalency}
\begin{split}
<u,v>_{\dot{H}^s(\mathbb{R}^n)}&=\int_{\mathbb{R}^n} \int_{\mathbb{R}^n} \frac{\displaystyle (\nabla u(x)-\nabla u(z)) \cdot (\nabla v(x)-\nabla v(z))}{\displaystyle |x-z|^{n-2+2s}} dx dz,\\
&=k(n,s)\int_{\R^n}|\xi|^{2s} \hat u(\xi)\overline{\hat v(\xi)}d\xi
\end{split}
\end{equation}
is an inner product on $\dot{H}^s(\mathbb{R}^n)$. Next, we introduce the notion of weak solution in $\dot{H}^s(\mathbb{R}^n)$ of the fractional Laplace equation.
\begin{definition}\label{frac_weak*}
For $1<s<2$ and  $u \in \dot{H}^s(\mathbb{R}^n)$, we introduce the weak formulation of  $(-\Delta)^s u$ as
\begin{equation}\label{frac_weak}
\langle (-\Delta)^s v,\phi \rangle=\frac{1}{k(n,s)}\langle v,\phi\rangle_{\dot H^s(\R^n)}\quad\text{for any }\phi\in C^\infty_0(\mathbb{R}^n).
\end{equation}
Here, $k(n,s)$ is as in  (\ref{inner_equivalency}).
\end{definition}

\section{The extension problem for the higher order fractional Laplacian}\label{prelim}
The well-known Caffarelli-Silvestre characterization of the fractional Laplacian $(-\Delta)^s$ for $0<s<1$ as a Dirichlet-to-Neumann map (established in \cite{CS07}) was generalized in \cite{Y13}, and more recently in \cite{CM22}, to all positive, non-integer orders of the fractional Laplace operator.  More specifically, for $1<s<2$, let  $b=3-2s$ and let $\Delta_b$ be the $b$-Laplace operator defined by
\begin{equation}\label{b-Laplace}
\Delta_b U=y^{-b} \nabla \cdot (y^b \nabla U).
\end{equation}
For $\Omega$ an open domain in $\mathbb{R}^{n+1}$, let $H^2 (\Omega,y^b)$ be the weighted Sobolev space equipped with the norm
 \begin{equation}\label{weight_norm}
 \|U\|^2_{H^2 (\Omega,y^b)}=\|y^{\frac{b}{2}} \Delta_b U\|^2_{L^2 (\Omega)}+\|y^{\frac{b}{2}} \nabla U\|^2_{L^2 (\Omega)}+\|y^{\frac{b}{2}}  U\|^2_{L^2 (\Omega)}.
 \end{equation}
Yang proves in \cite{Y13} the following result.
\begin{theorem}\textup{(\cite[Theorem 3.1]{Y13}, see also \cite[Theorem 1.2]{CM22})}\label{extension*}
Let $1<s<2$ and let $b=3-2s$.
Suppose $U$ is a function in $H^{2} (\mathbb{R}_{+}^{n+1},y^b)$ satisfying
\begin{equation*}
\left\{
\begin{array}{ll}
 \Delta^2_b U(x,y)=0 & \text { in } \mathbb{R}^n \times \mathbb{R}_{+}\\
U(x,0)=u(x) & \text { for all } x \in \mathbb{R}^n\\
\lim_{y \to 0} y^b U_{y}(x,y)=0 & \text { for all } x \in \mathbb{R}^n\\
\end{array}
\right.
\end{equation*}
for some $u\in \dot{H}^{s}(\mathbb{R}^n)$.
Here the operator $\Delta_b$ is as defined in (\ref{b-Laplace}), the PDE is satisfied in some weak sense, and the Dirichlet boundary condition is satisfied in the sense of traces.
Then,
\begin{equation}\label{D_to_N}
(-\Delta)^s u(x)=C_{n,s} \lim_{y \to 0} y^b\frac{\partial}{\partial y} \Delta_b U(x,y).
\end{equation}
\end{theorem}

Thanks to the Extension Theorem \ref{extension*}, we infer that boundary obstacle problems in $\mathbb{R}^n\times \{y\geq 0\}$ associated with the weighted bi-Laplace operator $\Delta^2_b$ and with an obstacle living on $\mathbb{R}^n \times \{0\}$ are equivalent to  obstacle problems for the fractional Laplacian $(-\Delta)^s$ for $1<s<2$.

\section{The obstacle problem for a higher order fractional Laplacian}\label{paper1}
In this section, we present an outline of our paper \cite{DHP}, joint with A. Petrosyan, where we study the regularity of solutions to  obstacle problems associated with the fractional Laplacian $(-\Delta )^s$ for $1<s<2$, in the Euclidean space $\mathbb{R}^n$ with $n \geq 2$. The problem consists in finding a function $u$ constrained to stay above a given obstacle $\psi$, and satisfying $(-\Delta)^s u=0$ in the region where it is above the obstacle. We assume that the obstacle $\psi$ is a $C^{1,1}$- function with compact support in the region under consideration.
In \cite{DHP}, we consider two versions of this  problem. In both cases, our solutions will live in the space $\dot{H}^s(\mathbb{R}^n)$ introduced in Section \ref{prelim}.
For the first version, we seek a function $u$, solution of the problem
\begin{equation}\label{prb1}
\left\{\begin{array}{ll}
u \geq \psi &\text{ in } \mathbb{R}^n\\
(-\Delta)^s u \geq 0 &  \text{ in }\mathbb{R}^n\\
(-\Delta)^s u =0  & \text{ in } \{u>\psi\}.\\
\end{array}\right.
\end{equation}
To analyze this problem, we follow a variational approach by minimizing the functional
\begin{equation}\label{prb1_fnct}
I_0[v]:=\int_{\mathbb{R}^n} \int_{\mathbb{R}^n} \frac{ \displaystyle |\nabla v(x)-\nabla v(z)|^2}{ \displaystyle |x-z|^{n-2+2s}} dx dz
\end{equation}
over the space of functions
\begin{equation}\label{prb1_admissible}
\mathcal{B}_0:=\{v \in \dot{H}^{s}(\mathbb{R}^n), \   v\geq \psi\}.
\end{equation}
We prove the existence and uniqueness of a minimizer $u_0$  of
(\ref{prb1_fnct}), then we show
that $u_0$ satisfies the variational inequality

\begin{equation*}
(-\Delta )^s u_0 \geq 0 \quad\mbox{in }\mathbb{R}^n
\end{equation*}
 in the weak sense
\begin{equation}\label{prb1_varineq}
\langle (-\Delta)^s u_0, v-u_0\rangle=\frac1{k(n,s)}\langle u_0, v-u_0\rangle_{\dot H^s}\geq 0\quad\text{for all }v \in \mathcal{B}_0.
\end{equation}

The second version is a problem with boundary values imposed on the complement of the unit ball in $\mathbb{R}^n$:
\begin{equation}\label{prb2_original}
\left\{\begin{array}{ll}
w \geq \Psi &\text{ in } \mathbb{R}^n\\
(-\Delta)^s w \geq 0 &  \text{ in } \mathbb{R}^n \\
(-\Delta)^s w =0  & \text{ in } B_1^{'} \cap \{w>\Psi\}\\
w(x)=g(x), \nabla w=\nabla g &\text { on } \mathbb{R}^n \setminus B^{'}_1. \\
\end{array}\right.
\end{equation}
Here $B^{'}_1$ is the unit ball in $\mathbb{R}^n$, $\Psi$ denotes the obstacle,  and  $g \in \dot{H}^{1+s}(\mathbb{R}^n)$ is a function satisfying $(-\Delta)^s g \geq 0$ and $g \geq \Psi$ in $\mathbb{R}^n$. Letting $u=w-g$, $f=-(-\Delta)^s g$, and $\psi=\Psi-g$, we see that (\ref{prb2_original}) is equivalent to
\begin{equation}\label{prb2}
\left\{\begin{array}{ll}
u \geq \psi &\text{ in } \mathbb{R}^n\\
(-\Delta)^s u \geq f(x)&  \text{ in } B^{'}_1\\
(-\Delta)^s u =f(x)  & \text{ in } B^{'}_1 \cap \{u>
\psi\}\\
u=0, \nabla u=0 &\text { on } \mathbb{R}^n \setminus B^{'}_1. \\
\end{array}\right.
\end{equation}
To construct a solution of (\ref{prb2}) we prove the existence and uniqueness of a minimizer $u$ of the functional
\begin{equation}\label{prb2_fnct}
I[v]:=\frac{1}{2 k(n,s)}\int_{\mathbb{R}^n} \int_{\mathbb{R}^n} \frac{\displaystyle |\nabla v(x)-\nabla v(z)|^2}{\displaystyle |x-z|^{n-2+2s}} dx dz - \int_{\mathbb{R}^n} f(x) \  (v(x)-\psi(x)) dx
\end{equation}
over the set
\begin{equation}\label{prb2_admissible}
\mathcal{B}:=\{v \in \dot{H}^{s}(\mathbb{R}^n),\   v=0, \nabla v=0 \text { on } \mathbb{R}^n \setminus B'_1,   \text { and } v\geq \psi\}.
\end{equation}
Here, the constant $k(n,s)$ is the same as in (\ref{inner_equivalency}). Next, we show that the minimizer $u$ satisfies

 \begin{equation*}
(-\Delta )^s u -f \geq 0\ \mbox{ in }\mathbb{R}^n
\end{equation*}
 in the weak sense
\begin{equation}\label{prb2_varineq}
\langle (-\Delta)^s u-f, v-u\rangle=\frac1{k(n,s)}\langle u, v-u\rangle_{\dot H^s}-\langle f, v-u\rangle_{L^2}\geq 0\quad\text{for all }v \in \mathcal{B}.
\end{equation}
In particular,  $\mu:=(-\Delta)^su-f$ is a measure on $B'_1$. To prove the regularity of our minimizer, we follow an approach from potential theory inspired by the one of Schild in \cite{Sc84} and \cite{Sc86}, as described in Section \ref{thin_section}. More specifically, we represent our minimizer and its Laplacian in terms of  Riesz potentials of some measures. Then, some of the main tools to reach the regularity results will be  a continuity principle and a maximum principle.

Applying the Extension Theorem \ref{extension*}, we see that  problems (\ref{prb1}) and (\ref{prb2}) are equivalent to  boundary obstacle problems in $\mathbb{R}^n\times \{y\geq 0\}$, with an obstacle living on $\mathbb{R}^n \times \{0\}$. In the case where $b=0$, the operator $\Delta_b^2$ reduces the bi-Laplacian, and the resulting problem is the object of Schild's work. Therefore, our paper can also be seen as a generalization of Schild's results to the case of the operator $\Delta^2_b$ for $ b \in (0,1)$.  For problem \eqref{prb1}, however, we are able to obtain an a.e.  representation of $u_0$ as the Riesz potential of order $2s$ of the measure $\mu_0$ on the whole space without resorting to the extension procedure. To illustrate our arguments,  we start by introducing  the fundamental solution $\phi_s$ of the fractional Laplacian $(-\Delta)^s$ given by (for $ x \in \mathbb{R}^n \setminus \{0\} $)
\begin{equation}\label{fundamental_fractional}
 \left\{\begin{array}{lll}
\phi_s(x)=\frac{1}{|x|^{n-2s}},&\   \text { when $n \neq 1,2,3$ or $s \neq \frac{3}{2}$},\\
\phi_{\frac{3}{2}}(x)=-\log|x|,&\  \text { when $n=3$,}\\
\phi_{\frac{3}{2}}(x)=-|x|,&\ \text { when $n=2$,}\\
\phi_{\frac{3}{2}}(x)=|x|^2\left(\log|x|-\frac{3+n}{2n+2}\right),&\ \text { when $n=1$.}
\end{array}\right.
\end{equation}

\begin{lemma}\label{u0_rep*}\textup{(\cite[Lemma 4.1]{DHP})}
The identity
\begin{equation}\label{u0_rep}
u_0=\phi_s*\mu_0
\end{equation}
 holds true for a.e. $x$ in $\mathbb{R}^n$.
\end{lemma}

For problem  (\ref{prb2}), instead, our representations in terms of  Riesz potentials of measures will be local.
 As a starting point, we extend our minimizer $u$ to the upper half space with a function $U\in H^2(\mathbb{R}^n,y^b)$  satisfying
\begin{equation}\label{extension_prb}
\left\{\begin{array}{ll}
\Delta_b^2 U=0 & \text { for } (x,y) \in \mathbb{R}^{n+1}_+\\
\lim_{y \to 0^+} y^b U_y (x,y)=0 &\text { for } (x,y) \in \mathbb{R}^{n+1}_+\\
U(x,0)=u(x)  & \text { for }x \in \mathbb{R}^n,\\
\end{array}\right.
\end{equation}
for $b=3-2s$.
Given the extension result of Theorem \ref{extension*}, we know that the extension function $U(x,y)$  of $u(x)$ solves
\begin{equation}\label{prb2_extension}
\left\{
\begin{array}{ll}
\Delta^2_b U =0 & \text { for } (x,y) \in B_1^{+}\\
U=0 , \  U_{\nu}=0 & \text { for } (x,y) \in (\partial B_1)^+ \\
\lim_{y \to 0} y^b U_{y}(x,y)=0  & \text { for all } x \in B_1'\\
U(x,0) \geq \psi(x) & \text { in } B_1'\\
\lim_{ y \to 0^+} y^b \frac{\partial}{\partial y} \Delta_b U(x,y) \geq f(x)  & \text{ in } B'_1\\
\lim_{ y \to 0^+} y^b \frac{\partial}{\partial y} \Delta_b U(x,y) =f(x)  & \text{ in } B'_1 \cap \{u>\psi\},\\
\end{array}
\right.
\end{equation}
Here  $B_1$ is the unit ball in $\mathbb{R}^{n+1}$, $B_1^+=B_1 \cap \{y>0\}$, $(\partial B_1)^+=\partial B_1 \cap \{y>0\}$, and $\nu$ is the outer unit normal to the ball $B_1$.
For fixed $0<\rho\leq 1/2$, let $\eta \in C_0^{\infty}(B'_{2 \rho})$ be such that $\eta \equiv 1$ on $B'_{\rho}$ and $|\nabla \eta| \leq \frac{1}{\rho}$ on $B'_{2\rho}$. We will denote by $\mu_{\rho}$ the restriction of the measure $\mu$ on $B'_{2\rho}$, that is

\begin{equation}\label{measure_ball}
\mu_{\rho}:=\eta \cdot \mu.
\end{equation}
We first obtain the following local representation of $u$ and $\Delta u$.

\begin{lemma}\label{rep_local*}\textup{(\cite[Lemma 6.2]{DHP})}
For a ball $B'_{\rho}$ in $\mathbb{R}^n$, there exists a weak solution $R_{\rho}(x)$ of $(-\Delta)^s R_{\rho}(x)= f$ in $B'_{\rho}$ such that
\begin{align}\label{u_rep_local}
&\bullet \quad u(x)=\phi_{s}*\mu_{\rho}(x)+ R_{\rho}(x) \quad\text{ \ for a.e.  $x \in B'_{\rho}$. }\\
&\bullet \quad \Delta u(x)=\Delta \phi_{s}*\mu_{\rho}(x)+ \Delta R_{\rho}(x)\quad\text{ for a.e. $x \in B'_{\rho}$. }\label{lapxu_rep_local}
\end{align}
Moreover, $R_{\rho}(x)$ satisfies the estimate
\begin{equation}\label{R_bnd}
\begin{split}
\|R_{\rho}\|_{W^{2,\infty}(B'_{\rho})} \leq C(n,s,\rho) &\left(\|y^{{b}/{2}} U\|_{L^2(B_{2 \rho})}  +\| y^{{b}/{2}} \nabla U\|_{L^2(B_{2 \rho})}\right.\\
&\left.+\|f\|_{L^2(B_{2 \rho})} \right).
\end{split}
\end{equation}
\end{lemma}
Denoting by $E_s$ the fundamental solution for the $b$-bi-Laplace equation defined by
\begin{equation}\label{fundamental}
 \left\{\begin{array}{ll}
E_s(x,y)=\frac{1}{\left(|x|^2+|y|^2\right)^{\frac{n-2s}{2}}},\qquad\qquad\qquad\ \text{when }n\neq 1,2,3,&\text{or }s \neq {3}/{2},\\
E_{{3}/{2}}(x,y)=-\log\left(|x|^2+|y|^2\right)^{\frac{1}{2}},&\text{when $n=3$,}\\
E_{{3}/{2}}(x,y)=-\left(|x|^2+|y|^2\right)^{\frac{1}{2}},&\text{when $n=2$,}\\
E_{{3}/{2}}(x,y)=\left(|x|^2+|y|^2\right)\left(\log\left(|x|^2+|y|^2\right)^{\frac{1}{2}}-\frac{3+n}{2n+2}\right),&\text{when $n=1$.}\\
 \end{array}\right.
\end{equation}
for $(x,y) \in \mathbb{R}^n \setminus \{(0,0)\}$.
The representations obtained in Lemma \ref{rep_local*} can be extended to the upper half space as follows.
\begin{lemma}\label{rep_local_upper*}\textup{(\cite[Lemma 9.1]{DHP})}
For every $y>0$, and  every ball $B'_{\rho}$ in  $\mathbb{R}^n$, the harmonic extension $U$ of $u$ and its Laplacian admit the following representations:

\begin{itemize}
\item For  every $y>0$ and almost every $x\in  B'_{\rho}$
\begin{equation}\label{u_rep_local_upper}
U(x,y)=\int_{\mathbb{R}^n} E_s(x-\xi,y) \ d\mu_{\rho}(\xi)+L_{\rho}(x,y),
\end{equation}
where $L_\rho(x,y)$ is the b-biharmonic extension of $R_\rho(x)$ to the upper half-space.
\item For  every $y>0$ and almost every $x\in  B'_{\rho}$
\begin{equation}\label{lapu_rep_local_upper}
\Delta U(x,y)=\int_{\mathbb{R}^n} \Delta E_s(x-\xi,y) \ d\mu_{\rho}(\xi)+\Delta L_{\rho}(x,y).
\end{equation}
\item For almost every $x\in  B'_{\rho}$
\begin{equation}\label{uyy_rep_local}
U_{yy}(x,0)= C(n,s,\rho) \Delta \phi _s *\mu_{\rho}(x) + (L_\rho)_{yy}(x,0).
\end{equation}
\item $L_{\rho}$ satisfies the estimate
\begin{equation}\label{L_bnd}
\begin{split}
\|L_{\rho}\|_{W^{2,\infty}(B_{\rho})}\leq C(n,s,\rho)& \left(\|y^{{b}/{2}} U\|_{L^2(B_{2 \rho})}\right.\\
& \left.  +\| y^{{b}/{2}} \nabla U\|_{L^2(B_{2 \rho})}+\|f\|_{L^{\infty}(B_{2 \rho})} \right).
\end{split}
\end{equation}
\end{itemize}
\end{lemma}
One of the main regularity results in our paper is the following $C^{1,1}_{ \rm loc}$-estimate of the extension function $U$ of the solution $u$.

\begin{theorem}\label{C11_reg*}\textup{(\cite[Corollary 7.9 \& Theorem 9.2]{DHP})}
The extension function $U$ is in $C^{1,1}_{\rm loc}(\overline{\mathbb{R}^{n+1}})$. In particular,
\begin{equation} \label{C1alpha_reg}
\begin{split}
\|U\|_{C^{1,1}(B^+_{{\rho}/{4}} \cup B'_{{\rho}/{4}})}\leq  C(n,s,\rho) &\left(\|\Delta \psi\|_{L^{\infty}(\mathbb{R}^n)}+\|y^{b/2} U\|_{L^2(B_{2 \rho})}\right.\\
&\left.+\| y^{b/2} \nabla U\|_{L^2(B_{2 \rho})}+\|f\|_{L^{\infty}(B_{2 \rho})} \right).\\
\end{split}
\end{equation}
\end{theorem}
The proof is based on several applications of a continuity principle and a maximum principle from potential theory.
 More specifically,  since $u$ is bounded on the support of the measure $\mu$, with the representation (\ref{u_rep_local}) in hand we are able to apply a continuity result from potential theory to get a bound for $u$ on every ball $B'_{\rho} \subset B'_1$. Next, we aim to prove the local boundedness of $\Delta_x u$. To this end, we first show its local boundedness on the support of the measure $\mu$ using lower semi-continuity properties and an application of Green's formula. Then, we adapt a maximum principle from potential theory to obtain local bounds for $\Delta_x u$ in the entire ball $B'_1$. Moreover, using some crucial properties of the mixed second order partial derivative of $\phi_s$, we are able to obtain some local representations of the ones of $u$, from which we infer the relevant local bounds. Finally, using  the representations in Lemma \ref{rep_local_upper*}, we are able to deduce the same $C^{1,1}_{loc}$-regularity estimate for the extension function $U$.
 
The second main regularity result in our paper is the following.

\begin{theorem}\label{H1+s_reg*}\textup{(\cite[Theorem 8.1]{DHP})}
Let $u$ be the minimizer of $I[\cdot]$. Then, $u \in H_{\rm loc}^{1+s}(\mathbb{R}^n)$.  In particular,
\begin{equation}
\label{H1+s_reg}
\begin{split}
\|u\|_{H^{1+s}(B'_{{\rho}/{8}})}\leq C(n,s,\rho) &\left(\|\Delta \psi\|_{L^{\infty}(\R^n)}+\|y^{{b}/{2}} U\|_{L^{2}(B_{2\rho})}\right.\\
& \left. +\|y^{{b}/{2}} \nabla U\|_{L^{2}(B_{2\rho})}+ \|f\|_{L^{\infty}(B_{2\rho})}\right).
\end{split}
\end{equation}
\end{theorem}
The proof hinges on the regularity result of Theorem \ref{C11_reg*}, and a  classical result from potential theory which gives a relation between the $L^2$-norm of the Riesz potential of order $\beta/2$ of a compactly supported measure $\nu$ and the Riesz energy of order $\beta$ of the same measure.
From Theorem \ref{H1+s_reg*}, Theorem \ref{C11_reg*}, and the representations of Lemma \ref{rep_local_upper*}, we infer the following $H^3_{loc}(\mathbb{R}^{n+1},y^b)$ regularity of the extension function $U$.
\begin{corollary}\label{H3_reg*}\textup{(\cite[Corollary 9.4]{DHP})}
Let $u$ be the minimizer of $I[\cdot]$  and let $U$ be the extension
of $u$ to the upper half space $\mathbb{R}_+^{n+1}$. Then, $U\in
H_{\rm loc}^{3}(\mathbb{R}^{n+1},y^b)$. Moreover,
\begin{equation}\label{H1+s_to_H3}
\begin{split}
\|U\|_{H^3(B_{\rho/2}, y^b)}\leq &C(n,s,\rho)) \left(\|\Delta \psi\|_{L^{\infty}(\R^n)}+\|y^{{b}/{2}} U\|_{L^{2}(B_{16\rho})}\right.\\
&\left.+\|y^{{b}/{2}} \nabla U\|_{L^{2}(B_{16\rho})} +\|y^{{b}/{2}} \Delta_b U\|_{L^{2}(B_{\rho/2})}+ \|f\|_{L^{\infty}(B_{16\rho})}\right).
\end{split}
\end{equation}
\end{corollary}

\section{A two-phase boundary obstacle-type problem for the bi-Laplacian}

In this section, we present an outline of our recent work \cite{DH22}, where we consider the following two-phase boundary obstacle problem
\begin{equation}\label{prb_2phase}
\left\{
\begin{array}{ll}
\Delta^2 u =0 & \text { in } B_1^+\\
u=g & \text { on } (\partial B_1)^+ \\
u_y=0 & \text { on } B_1'\\
(\Delta u)_{y}=\lambda_- (u^-)^{p-1}-\lambda_+ (u^+)^{p-1}  & \text { on }  B_1'.\\
\end{array}
\right.
\end{equation}
Here  $B_1^+ =\{z=(x,y) \in B_1 ; \ \ y>0\}$, with $B_1$ being the unit ball in $\mathbb{R}^{n+1}$ centered at the origin, $ (\partial B_1)^+= \partial B_1 \cap \{y > 0\}$, $B_1'=B_1 \cap \{y=0\}$, $\lambda_-, \lambda_+$ are positive constants, and $p>1$. For $q:=\max\{2,p\}$, we let $g\in W^{2,q}(B_1^{+})$, with $g_y=0$ on $B_1'$.

The  problem (\ref{prb_2phase}) has been considered in \cite{Y13} and \cite{FF20} in the case where the non-homogeneous boundary condition is given by $(\Delta u)_y=0$ in \cite{Y13} and $(\Delta u)_y= h(x) u(\cdot,0)$ in \cite{FF20}, with $h(x)$ being a $C^1$-function on $B_1'$. In both papers, the authors prove a monotonicity formula of Almgren's type and classify the possible blow-up limit profiles. Then, they use their results to establish the strong unique continuation property for a solution.

In \cite{DH22}, both the formulation of  problem (\ref{prb_2phase}) and our objectives are different. The non-homogeneous thin boundary condition $(\Delta u)_{y}=\lambda_- (u^-)^{p-1}-\lambda_+ (u^+)^{p-1}   \text { on }  B_1'$ has a right-hand-side which is non-differentiable in $u$ when $1<p\leq 2$ and non-linear in $u$ when $p>2$. We establish the optimal regularity of the solution of (\ref{prb_2phase}) for all $p>1$.  Moreover, in the cases when $p=2$ and $p\geq 3$, we study the structural properties of the free boundary $\left(\partial \{u>0\} \cup \partial \{u<0\}\right) \cap B_1'$.

To construct our solution, we minimize the functional
\begin{equation}\label{fnct_2phase}
J[w]:=\int_{B_1^+}  (\Delta w)^2 dx dy +  \frac{2}{p} \int_{B_1'} \lambda_- (w^-)^p+\lambda_+ (w^+)^p dx
\end{equation}
over the admissible set
\begin{equation}\label{admissible_2phase}
\mathcal{A}:= \left\{ w \in W^{2,q} (B_1^+), \ w=g \text { on } (\partial B_1)^+, \text { and } \ w_y=0 \text { on } B_1'\right\}.
\end{equation}
We show that the unique minimizer $u$ of (\ref{fnct_2phase}) over (\ref{admissible_2phase}) is a weak solution of (\ref{prb_2phase}) in the sense
\begin{equation}\label{weaksol}
 \int_{B_1^+} \Delta u \Delta \phi dx dy =\int_{B_1'}\left(\lambda_- (u^-)^{p-1}  -\lambda_+ (u^+)^{p-1}\right)\phi \ dx
 \end{equation}
  for all functions $\phi \in C^{\infty} (B_1^+)$ such that  $\phi=0 \text { on }  (\partial B_1)^+$ and $\phi_y=0 \text { on }  B_1'$.
  \subsection{Regularity of the minimizer}\label{reg_2phase_section}
The first step in studying the regularity of the  minimizer $u$ is to  establish the weak differentiability of $\Delta u$.  The proof is based on an application of the Lax-Milgram Theorem.
\begin{theorem}\label{gradv}\textup{(\cite[Theorem 2.5]{DH22})}
Let $u$ be the minimizer of (\ref{fnct_2phase}). Then $\Delta u$  is weakly differentiable, and its weak derivatives are in $L(B^+_{\rho})$ for all balls $B_{\rho} \subset \subset B_1$.
\end{theorem}
We then establish the following optimal regularity result.
\begin{theorem}\label{reg_2phase}\textup{(\cite[Lemma 3.3 \& Theorem 3.4]{DH22})}
Let $u$ be the minimizer of (\ref{fnct_2phase}). When $p$ is an integer,
 $u \in C_{loc}^{p+1, \alpha}(B_1^+\cup B_1')$ for all $\alpha<1$. Moreover,
 \begin{equation}\label{u_part_reg_norm1}
\|u\|_{C^{p+1,\alpha}(B_{\rho}^+\cup B_{\rho}^{'}) } \leq C(\alpha,\rho) \left(\|u\|_{L^{\infty}(B_{\rho})}+\|\Delta u\|_{L^{\infty}(B_{\rho})}\right).
\end{equation}
When $p$ is not an integer, $u \in C_{loc}^{\lfloor p+1 \rfloor, \alpha}(B_1^+)\cup B_1'$  for all $\alpha<p-1-\lfloor p-1 \rfloor$ and
\begin{equation}\label{u_part_reg_norm2}
\|u\|_{C^{p+1,\alpha}(B_{\rho}^+\cup B_{\rho}^{'})} \leq C(\alpha,\rho) \left(\|u\|_{L^{\infty}(B_{\rho})}+\|\Delta u\|_{L^{\infty}(B_{\rho})}\right).
\end{equation}
Furthermore, in the case where $p$ is an integer, if $z_0$ is a free boundary point such that $\nabla_x u (z_0) \neq 0$, then $u$ is not $C^{p+1,1}$ at $z_0$.
\end{theorem}
To prove Theorem \ref{reg_2phase}, we start with proving the local boundedness of the solution and its Laplacian using the subharmonicity of their positive and negative parts. Then, we achieve the $C_{loc} ^{1,\alpha}$-regularity of the solution ($0<\alpha<1$)  using {Morrey's method}. In particular, $u$ is locally Lipschitz in $B_1$, and hence
\begin{itemize}
		\item For  $1<p < 2$,  the function $\lambda_- (u^-)^{p-1}-\lambda_+ (u^+)^{p-1}$ is Holder continous of order $\alpha$ for all $\alpha<p-1$,
		\item For $p \geq 2$, and $p$ is an integer,  the function $\lambda_- (u^-)^{p-1}-\lambda_+ (u^+)^{p-1}$ is in $C_{loc}^{p-2, \alpha}(B_1)$ for all $0<\alpha<1$,
		\item For $p \geq 2$, and $p$ is not an integer,  $u \in C_{loc}^{\lfloor p-2 \rfloor, \alpha}(B_1)$ for all $\alpha<p-1-\lfloor p-1 \rfloor$.\\
\end{itemize}
At this point, we apply the regularity theory of oblique derivative problems to the Laplacian $\Delta u$ of our solution to conclude that $\Delta u \in C_{loc}^{p-1, \alpha}(B_1)$ for all $\alpha<1$ when $p$ an integer, and  $\Delta u \in C_{loc}^{\lfloor p-1 \rfloor, \alpha}(B_1^+)$  for all $\alpha<p-1-\lfloor p-1 \rfloor$ when $p$ is not an integer. Applying the same regularity theory, but this time to the solution $u$, we reach the desired regularity. Finally, we prove that the regularity achieved in (\ref{reg_2phase}) is optimal by arguing by contradiction.

 \subsection{Regularity of the free boundary}Here and in the sequel, $v$ will denote the Laplacian of the solution $u$, that is, $v:=\Delta u$.
To study the regularity of the free boundary
\begin{equation}\label{free_bndry}
\begin{split}
&\Gamma(u):=\Gamma_-(u) \cup \Gamma_+(u), \\
&\text{ with }
\Gamma_-(u)=\partial_x (\{u<0\} \cap B_1')
\text{ and }\Gamma_+(u)=\partial_x (\{u>0\} \cap B_1')\\
\end{split}
\end{equation}
we define the \emph{regular set} of the free boundary as
\begin{equation}\label{regularset}
\mathcal{R}(u)=\{z=(x,0) \in \Gamma(u) \ \ | \  u(x,0)=0, \nabla_x u (x,0) \neq 0,  \nabla_x v (x,0) \neq 0 \},
\end{equation}
and the\emph{ singular set}  $\Sigma(u):=\Gamma (u) \setminus \mathcal{R}(u)$.
The following regularity of the free boundary around regular points is a straightforward consequence of the regularity of the solution and of the implicit function theorem.
\begin{theorem}\label{regular_fb}\textup{(\cite[Theorem 6.1]{DH22})}
If $(x,0) \in \mathcal{R}(u)$, then the free boundary is a $C^{3,\alpha}-$ graph in a neighborhood of $(x,0)$  for all $0<\alpha<1$.
\end{theorem}
To study the behavior of the free boundary at singular points, we need to carry out blow-up analysis adapting approaches developed in analyzing the free boundary in the \emph {Signorini problem} (see for instance  \cite{ACS08}, \cite{CSS08}, \cite{CS07} \cite{GP09}, the book \cite{PSU12}), and the survey \cite{DS17}.

\subsection{Almgren's monotonicity, growth estimate and blow-up analysis at singular free boundary points}
We start with introducing \emph{Almgren's Frequency Functional} for functions $u, v\in W^{1,2}_{loc}(B_1)$:
\begin{equation}\label{Almgren's}
N_0(r,u,v):=r \frac{\displaystyle \int_{B_r^+} (|\nabla u|^2+|\nabla v|^2) dx dy}{\displaystyle \int_{(\partial B_r)^+} (u^2+ v^2) dS}.
\end{equation}
Then, for $u$ the solution of our problem (\ref{prb_2phase}), we aim to show that $N_0(r):=N_0( r, u, \Delta u)$ has a non-negative limit as $r \to 0^+$. To this end, we first prove that  a suitable perturbation of $N_0(r)$  has a limit  as $r \to 0^+$. Then we are able to show that $N_{0}(r)$ also has a limit as $r \to 0^+$, which coincides with  $\lim_{r \to 0^+} N(r)$. More precisely,  we prove the following result.

\begin{theorem}\label{Almgren's*} \textup{(\cite[Theorem 4.1 \& Corollary 4.4]{DH22})}
For $u, v\in W^{1,2}_{loc}(B_1)$,  define the \emph{perturbed Almgren's Frequency Functional}
\begin{equation}\label{Almgren's_per}
\begin{split}
N(r,u,v):=&r \frac{\int_{B_r^+}\displaystyle (|\nabla u|^2+|\nabla v|^2 + uv) dx dy}{\displaystyle \int_{(\partial B_r)^+} (u^2+ v^2) dS}+\\
&\frac{\displaystyle \int_{B_r{'}} \left(\lambda_- (u^-)^{p-1}  -\lambda_+ (u^+)^{p-1}\right)v dx}{\displaystyle \int_{(\partial B_r)^+} (u^2+ v^2) dS}.
\end{split}
\end{equation}
Also, for $u$ solution to (\ref{prb_2phase}) with $p=2$ or $p\geq 3$, let $N(r):=N( r, u, \Delta u)$. Then,
$$\mu:=\lim_{r \to 0^{+}} N(r) \text { is finite and } \mu \geq 0.$$ Moreover, the limit of $N_{0}(r)$ also exists and coincides with  $\lim_{r \to 0^+} N(r)$.
\end{theorem}

Next, we introduce the quantity
\begin{equation}\label{phi}
\phi(r) =\frac{1}{r^{n}}\int_{(\partial B_r)^+} \left(u^2+ v^2\right) dS,
\end{equation}
and we obtain the following estimates for it.

\begin{lemma} \label{phi_est*}\textup{(\cite[Lemma 4.5]{DH22})}
Let $\mu:=\lim_{r \to 0^+} N(r)$. Then, the following hold:
\begin{enumerate}
\item There exists a constant $\alpha>0$, and $r_0>0$ such that
\begin{equation}\label{phi_est_upper}
r^{-2\mu} \phi(r) \leq R^{-2\mu} \phi(R)  e^{\frac{2C(R^{\alpha}-r^{\alpha})}{\alpha}}
\end{equation}
 for all $0<r<R<r_0$.
\item For any $\delta>0$, there exists $r_0(\delta)$ such that

\begin{equation}\label{phi_est_lower}
r^{-2\mu} \phi(r) \geq R^{-2\mu} \phi(R) \left(\frac{r}{R}\right)^{2\delta}
\end{equation}
for all $0<r<R<r_0(\delta)$.
\end{enumerate}
\end{lemma}
As a corollary of Lemma \ref{phi_est*}, we obtain the following growth estimate.
\begin{corollary} \textup{(\cite[Corollary 4.6]{DH22} \label{growth})}
Let $u$ be the minimizer of (\ref{fnct_2phase}). Then, for all $z \in B^+_r\cup {B_r^{'}}$, $0<r<\frac{1}{2}$, we have
$$u(z) \leq C r^{\mu} \ \  \text { and } \ \ v(z) \leq C r^{\mu}, $$
where $C$ depends on the dimension and the local $L^{\infty}$-norms of $u$ and $v$.
\end{corollary}
The next step is to study blow-up sequences around a free boundary point $z_0=(x_0,0)$. Without loss of generality, we assume $z_0=0$ and let $\mu:=\lim_{r \to 0^+} N(r)$. For $0<r<1$ and $z=(x,y) \in B_1^+\cup B_1^{'}$, we define the \emph{Almgren's rescalings}
\begin{equation}\label{Almgren's_rescaling}
u_r(z):= \frac{u(rz)}{\sqrt{\phi(r)}} \text { and } v_r(z):= \frac{ v(rz)}{\sqrt{\phi(r)}}.
\end{equation}
Using the monotonicity result Theorem \ref{Almgren's*} and the estimates of the function $\phi$ established in Lemma \ref{phi_est*}, we prove the convergence, up to a subsequence, of Almgren's rescalings  $u_r$ and $v_r$ around singular points  to  homogenous harmonic polynomials $p_\mu$ and $q_\mu$, respectively,  of degree $\mu\geq 2$.
\begin{theorem}\label{blow-up*}\textup{(\cite[Theorem 5.1]{DH22})}
Let $u$ be the solution to (\ref{prb_2phase}) with $p=2$ or $p\geq 3$, and $v=\Delta u$. There exists subsequences $u_j:=u_{r_j}$ and $v_j:=v_{r_j}$ such that $u_j \to u^{\ast}$, $\nabla u_j \to \nabla u^{\ast}$, $v_j \to v^{\ast}$ and $\nabla v_j \to \nabla v^{\ast}$ as $j\to\infty$, uniformly on every compact subset of $\overline{\mathbb{R}_+^{n+1}}$. Moreover, $u^{\ast}$ and $v^{\ast}$ are homogenous harmonic polynomials of degree $\mu$. If, in addition, we have $\nabla_x u(0)=0$ or $\nabla_x v(0)=0$, then we must have $\mu \geq 2$.
\end{theorem}

\subsection{Monneau-type monotonicity formula and non-degeneracy of the solution at singular free boundary points}

To continue our analysis of the structure of the singular set, we need to derive some crucial proporties such as non-degeneracy, uniqueness of blow-ups, and  continuous dependence  of blow-ups on the singular set. A crucial tool in this endeavor is the monotonicity of a Monneau-type functional.
We define the Monneau-type functional
 \begin{equation}\label{Monneau}
M_{\mu} (r,u, v, p_{\mu}, q_{\mu})=\frac{1}{r^{n+2\mu}} \int_{(\partial B_r)^+} \left((u-p_{\mu})^2 + (v-q_{\mu})^2\right) dS,
\end{equation}
 where $p_{\mu}$ and $q_{\mu}$ are two homogeneous harmonic polynomials of degree $\mu$. To prove the monotonicity of $M_{\mu}$, we  introduce the Weiss-type functional
\begin{equation}\label{Weiss}
W_{\mu} (r,u,v)=\frac{H(r,u,v)}{r^{n+2\mu}} (N_0(r,u,v)-\mu).
\end{equation}
We have  the following result.

\begin{theorem}\label{Monneau*}\textup{(\cite[Theorem 6.2 \& Corollary 6.3]{DH22})}
Let $u$ be the solution to (\ref{prb_2phase}), and let $p_{\mu}$ and $q_{\mu}$ be two homogenous harmonic polynomials of degree $\mu$, symmetric with respect to the thin space $\{y=0\}$. Then, there exists positive constants $C_1$ and $C_2$ such that

\begin{equation}\label{monneauest1}
\frac{d}{dr} M_{\mu} (r,u, v, p_{\mu}, q_{\mu}) \geq \frac{2}{r}W_{\mu} (r,u,v)-C_1 \geq -C_2
\end{equation}
 In particular, $M_{\mu}(r,u, v, p_{\mu}, q_{\mu})$ has a limit as $r \to 0^+$.
\end{theorem}
Next, we use the monotonicity result of Theorem \ref{Monneau*}
to prove the non-degeneracy of the solution at singular free boundary points.
\begin{lemma}\label{non_degen*}\textup{(\cite[Lemma 6.4]{DH22})}
Let $u$ be the solution to (\ref{prb_2phase}) with $u(0)=0$ and either $\nabla _{x} u(0)=0$ or $\nabla _{x} v(0)=0$. There exists $C>0$ and $0<R_0<1$, depending possibly on $u$, such that
\begin{equation}\label{non_deg}
\sup_{(\partial B_r)^+} |u| \geq C r^{\mu}\text { or } \sup_{(\partial B_r)^+} |v| \geq C r^{\mu}
\end{equation}
for all $0<r<R_0$.
\end{lemma}

At this point, we are ready to study the dimension of the singular set.

\subsection{The structure of the singular set}
Let $z_0$ be a point on the free boundary of $u$, the solution of (\ref{prb_2phase}).
For $R\geq 1$, we introduce the \emph{homogenous rescalings}
\begin{equation}\label{homo_rescaling}
u_r^{\mu}(z):=\frac{u(z_0+rz)}{r^{\mu}} \text { and } v_r^{\mu}(z):=\frac{v(z_0+rz)}{r^{\mu}},
\end{equation}
for   $z \in B_R^+ \cup B^{'}_R\text { and } 0<r<\frac{1}{R}$.
Using the growth estimate of Corollary \ref{growth} and the non-degeneracy result of Lemma \ref{non_degen*}, we prove the convergence of these rescalings to two harmonic polynomials homogenous of degree $\mu$. Moreover, we prove the uniqueness of the blow-up polynomials using the monotonicity of Monneau's functional (Theorem \ref{Monneau*}). More precisely, we have the following.

\begin{theorem}\label{homo_blow-up}\textup{(\cite[Theorem 7.1]{DH22})}
Let $u$ be the solution of (\ref{prb_2phase}) and suppose that $z_0$ is a free boundary point such that $\nabla _{x}u(z_0)~=~0$ or $\nabla _{x}v(z_0)=0$. Then, there exist two unique harmonic polynomials $p_\mu$ and $q_\mu$, homogenous of degree $\mu$  and such that $u_r^{\mu} \to p_\mu$  in $W_{loc}^{3,2}(B^+_R \cup B_R^{'})$ and in $C^{3, \alpha}(B_R^+ \cup B_R^{'})$, and $v_r^{\mu} \to q_\mu$ in $W_{loc}^{1,2}(B^+_R)$ and in $C^{1, \alpha}(B_R^+\cup B_R^{'})$ for all $R\geq 1$.
\end{theorem}

For a free boundary point $z_0$, let \  $N_0^{z_0}$ denote Almgren's frequency formula at the free boundary point $z_0$,\  $\mu_{z_0}=\lim_{r \to 0} N_0^{z_0}$, and
$p_{\mu}^{z_0}$ and $q_{\mu}^{z_0}$ the $\mu-$homogenous blow-up polynomials at $z_0$ as in Theorem \ref{homo_blow-up}.
For $\mu \in \mathbb{N}$,  we define
\begin{equation}\label{mu_singular}
\begin{split}
\Sigma_{\mu}\{&(x,0) \in \Gamma(u) \ \  |\\
& \ \  u(x,0)=0,  \nabla_x u (x,0)=0 \text { \emph{or} } \nabla_x v (x,0)=0,  \text { and }   N_0^{z_0} (0+,u,v)=\mu\}.
\end{split}
\end{equation}
We denote by $\mathcal{P}(u)$ the class of homogenous harmonic polynomials of degree $\mu$ and even in $y$. We note that $\mathcal{P}(u)$ is a convex subset of the finite dimensional vector space of all polynomials homogenous of degree $\mu$. We endow  $\mathcal{P}(u)$ with the norm of  $L^2((\partial B_1)^+)$.

The monotonicity of Monneau functional yields the  continuous dependence of the blow-ups on the singular points as in Theorem \ref{homo_blow-up}.

\begin{theorem}\label{cont_dep}\textup{(\cite[Theorem 7.2]{DH22})}
Let $u$ be the solution to (\ref{prb_2phase}). The mapping $z_0 \mapsto (p_{\mu}^{z_0}, q_{\mu}^{z_0})$ from $\Sigma_{\mu}$ to $\mathcal{P}(u) \times \mathcal{P}(u)$  is continuous. Moreover, for any compact set $K \subset \subset \Sigma_{\mu}$, there exists a modulus of continuity $\omega_\mu$ such that
$$|u(z)-p^{z_0}_{\mu}(z-z_0)|\leq \omega_\mu (|z-z_0|) |z-z_0|^{\mu}$$
for any $z_0 \in K$.
\end{theorem}

To understand the dimension of the singular set, we first establish the following result. The proof is based on the growth estimate Corollary \ref{growth} and the non-degeneracy result  Lemma \ref{non_degen*}.

\begin{lemma}\label{count_union}\textup{(\cite[Theorem 7.3]{DH22})}
Let $u$ be the solution to (\ref{prb_2phase}). Then $\Sigma_{\mu}$ is of type $F_{\sigma}$, that is, it is a union of countably many closed sets.
\end{lemma}

At this point, we are ready to establish the rectifiability of the singular set. We denote by $d^{z_0}_{\mu}$ the dimension of $\Sigma_{\mu}$ at a point $z_0 \in \Sigma_{\mu}$. That is,
\begin{equation}\label{mu_sing_dim}
d^{z_0}_{\mu}=\max \{d^{u,z_0}_{\mu}, d^{v,z_0}_{\mu}\},
\end{equation}
where $$d^{u,z_0}_{\mu}=dim\{\zeta \in \mathbb{R}^{n} \ \ |\ \ \zeta \cdot \nabla_x p^{z_0}_{\mu}(x,0)=0 \text { \ \  for all x } \in \mathbb{R}^n\}$$
and
$$d^{v,z_0}_{\mu}=dim\{\zeta \in \mathbb{R}^{n} \ \ | \ \ \zeta \cdot \nabla_x q^{z_0}_{\mu}(x,0)=0 \text { \ \  for all x } \in \mathbb{R}^n\}.$$
We also let  $\Sigma^d_{\mu}:=\{ z_0 \in \Sigma_{\mu} (u) | d^{z_0}_{\mu}=d \}$.  The proof  of the following rectifiability result is based on the continuous dependence on blow-ups,  Whitney's extension theorem, and the implicit function theorem as in \cite[Theorem 1.3.8]{GP09}.
\begin{theorem}\label{singular_dim}\textup{(\cite[Theorem 7.4]{DH22})}
Let $u$ be the solution to (\ref{prb_2phase}). Then for every $\mu \in \mathbb{N}$ and $d=0,1,2,..,n-2$, the set $\Sigma^d_{\mu}$ is contained in a countable union of d-dimensional $C^1$-manifold.
\end{theorem}




\strictpagecheck
\checkoddpage
\ifoddpage
\newpage{\ }\thispagestyle{empty}
\fi

\end{document}